\documentclass{setup_class}

\title[Adaptive Constraint Sampling for SGD]{Adaptive Constraint Sampling for Stochastic Gradient Methods in Constrained Optimization}

\usepackage[short]{optidef}
\usepackage{bbm}
\usepackage{tikz}
\usetikzlibrary{arrows.meta, shapes.geometric, positioning, calc, shapes.misc}
\usepackage{algorithm}
\usepackage{subcaption}
\usepackage{pgfplots}
\usepackage{pgffor}
\pgfplotsset{compat=1.18}
\usepackage{hyperref}
\usepackage{enumitem}
\usepackage{cleveref}
\usepackage{booktabs}
\usepackage{adjustbox}
\usepackage{microtype}

\newtheorem{assumption}{Assumption}

\newcommand{\R}{\mathbb{R}}

\newcommand{\norm}[1]{\left\lVert #1 \right\rVert}
\newcommand{\G}{\mathbb{G}}

\newcommand{\Scal}{\mathcal{S}}
\newcommand{\Acal}{\mathcal{A}}
\newcommand{\Ucal}{\mathcal{U}}

\newcommand{\Gikzero}{\mathbb{G}_{I,0}}
\newcommand{\Gik}{\mathbb{G}_{I,k}}
\newcommand{\Gikk}{\mathbb{G}_{I,k+1}}

\newcommand{\N}{\mathbb{N}}
\newcommand{\m}{\overline{1:m}}

\newcommand{\Dcal}{\mathcal{D}}

\newcommand{\Kcal}{\mathcal{K}}

\newcommand{\jk}{{j_k}}

\newcommand{\n}{\overline{1:n}}

\author{%
 \Name{Erik M\"uller} \Email{erik.peter.mueller@rwth-aachen.de}\\
 \addr Chair of Intelligent Control Systems, RWTH Aachen University, 52062, Aachen, Germany
 \AND
 \Name{Naum Dimitrieski} \Email{naum.dimitrieski@ic.rwth-aachen.de}\\
 \addr Chair of Intelligent Control Systems, RWTH Aachen University, 52062, Aachen, Germany%
 \AND
 \Name{Mohamed-Ali Belabbas} \Email{belabbas@illinois.edu}\\
 \addr Coordinated Science Laboratory,
University of Illinois at Urbana-Champaign, Urbana, IL 61801~USA%
 \AND
 \Name{Christian Ebenbauer} \Email{christian.ebenbauer@ic.rwth-aachen.de}\\
 \addr Chair of Intelligent Control Systems, RWTH Aachen University, 52062, Aachen, Germany%
}

\begin{document}

\maketitle

\begin{abstract}%
This paper proposes an adaptive constraint-sampling technique for stochastic gradient-based methods to efficiently solve constrained optimization problems with a large number of inequality constraints. 
The proposed constraint-sampling technique maintains a set of recently violated constraints and draws the constraint used in each iteration uniformly either from this set or from the set of all constraints. This keeps its per-iteration arithmetic overhead independent of the number of constraints. In simulations on synthetic convex optimization problems and in the synthesis of a control barrier function for the Dubins Car, the constraint-sampling technique is shown to significantly outperform uniform constraint sampling.
\end{abstract}


\section{Introduction}\label{sec:introduction}
Constrained optimization problems with many inequality constraints can be found in many fields of control and learning. Examples include learning-based control with synthesis of control barrier functions (CBFs) from existing data \citep{robey2020learning}, or robust and stochastic control, where each existing sampled realization of the (stochastic) system provides an additional constraint \citep{calafiore2006scenario, campi2009scenario}. A further example is approximate dynamic programming, with one constraint present for each sampled state-action pair \citep{de2003linear, de2004constraint, samiedaluie2024approximate}. In the field of machine learning, constrained formulations are also known to be used (see, e.g., \citet{kotary2021ijcai, gallego2025cooper}), with one particular application found in adversarial robustness (see, e.g., \citet{robey2021adversarial, rezkellah2025machine}).

One approach to solve such problems consists of so-called active-set or working-set methods, as discussed in, e.g., \citet{gould2005numerical} and the references therein. Such methods typically estimate the active set or maintain a working set of constraints, and then perform an optimization step or solve the reduced problem, in a sequential manner (see, e.g., \citet{tits2006constraint,schittkowski2009anactive,laiu2019aconstraint}). Further, identifying the active set has also been studied in the literature (see, e.g., \citet{burke1988on}, \citet{facchinei1998onthe}, \cite{duchi2021asymptotic} and the references therein). However, such methods typically require at least periodic evaluations of all constraints, which might be computationally expensive in large-scale problems. Recent works, such as \citet{li2025new}, try to tackle this problem and focus on developing more efficient methods. Similar problems have arisen in working-set methods for semi-infinite constrained problems. To that end, approaches have been proposed which update the working set with information of only a subset of all constraints (see, e.g., \citet{heunis1985use,wardi1988stochastic}). However, updating the working set requires also solving a subproblem and evaluating the existing constraints, which can become computationally expensive as well. Additionally, approaches that consider reduced problems based on randomly subsampled constraints have also been discussed in the literature on random convex programs (see, e.g., \citet{calafiore2010random}).

An alternative to sequentially solving subproblems consists of iterative optimization methods that randomly sample the constraints at every iteration. Such methods are developed and discussed in, e.g., \citet{polyak2001random,nedic2011random,chakraborty2026randomized}. Further, \citet{wang2016stochastic} extends this approach to stochastic first-order methods, with alternative approaches presented in \citet{xu2020primal,nedich2023huber,dimitrieski2025stochastic}. Typically, these methods use uniform constraint sampling.
However, in a large-scale problem, there may be many constraints which do not provide relevant information for the current update step. Introducing constraint subsampling, including uniform subsampling, to a first-order method can lead to very 'noisy' updates, specifically, from sampling currently 'irrelevant' constraints. This usually leads to worse algorithm performance, and motivates designing adaptive constraint-sampling techniques (CSTs) 
%in particular, ones 
whose implementation does not depend on the number of constraints, and which provide 'useful' constraint information. 

Adaptive constraint-sampling approaches have been developed in the convex optimization literature, for instance, in \citet{clarkson1988las,clarkson1995vegas}, \citet{wei2020arandomized}, and \citet{clarkson2012sublinear}, \citet{cotter2016light}.  The former two iteratively solve subproblems in order to adapt the constraint-sampling probabilities, which requires evaluating all the constraints. A further, dual-based approach based on solving subproblems is given in \citet{wei2020arandomized}. The latter two do not solve subproblems, rather, they are iterative first-order based methods which adapt the sampling probabilities for each gradient-based step, with an $\mathcal{O}(m)$ per-iteration arithmetic overhead. Finally, \citet{sonthalia2022project} proposes a Bregman-projection based method with adaptive constraint sampling, where a working set can be updated solely from randomly sampled constraints. However, it requires full objective information and Bregman projections, which may be computationally expensive in large-scale problems.

\textbf{Contributions:} $\mathbf{1)}$ We present a novel adaptive CST for stochastic gradient-based methods. The adaptation mechanism is working-set based. In contrast to classical working-set methods, it does not consider the full set of constraints for updating the working set, but uses a history-dependent mechanism based on the sampled constraints thus far. Moreover, it
does not require full objective information nor solving subproblems for every working-set update, and is thus compatible with stochastic gradient-based methods. Lastly, the CST can be implemented with a per-iteration arithmetic overhead of amortized~$\mathcal{O}(1)$.

$\mathbf{2)}$ We implement the proposed CST on four stochastic gradient-based methods adapted from \citet{wang2016stochastic, nedich2023huber, dimitrieski2025stochastic, xu2020primal}, and numerically show the advantage of using the CST over uniform constraint sampling on synthetic optimization problems and applications in learning CBFs for the Dubins Car.

\section{Problem Statement and Solution Approach}

In this paper, we consider solving a constrained optimization problem of form
\begin{equation}
    x^{*} \in \mathrm{argmin}_{x\in \R^d} \hspace{0.1cm}  f(x) \hspace{0.4cm} \mathrm{s.t.} \hspace{0.4cm} g_j(x) \leq 0, \hspace{0.4cm}  j \in \overline{1:m}, \label{problem:1}
\end{equation}
where~${\overline{1:m} \coloneqq \{1,2,\dots,m\}}$ and $m$ is very large. We consider the finite-sum objective to be given as $${f(x) := \tfrac{1}{n}\sum\nolimits_{i=1}^{n}f_i(x)}$$
with components ${f_i: \R^d \to \R}$ for all ${i \in \overline{1:n}}$, and Problem~\eqref{problem:1} to have affine constraints, i.e.,  ${g_j(x) \coloneqq a_j^\top x +b_j}$ with $a_j \in \R^d \setminus \{0\}$ and $b_j\in \R$ for all~${j \in \overline{1:m}}$. We define the active set~as
\begin{equation*}
    \Acal\coloneqq\{j \in \m : g_j(x^{*}) = 0 \},
\end{equation*}
where we assume $x^*\in \R^d$ to be the unique minimizer of Problem~\eqref{problem:1}.

To solve Problem~\eqref{problem:1}, i.e., find $x^* \in \mathbb \R^d$, we employ four stochastic gradient-based methods adapted from \cite{wang2016stochastic, nedich2023huber, dimitrieski2025stochastic, xu2020primal} with the following update steps:
\begin{subequations}
\begin{itemize}[]
    \item Projection-SGD \citep[Algorithm 1, Eq. (9) with $\bar{x}_k = x_k$]{wang2016stochastic}
    \begin{equation}\label{alg:projection}
    z_k = x_k - \alpha_k \nabla f_{i_k}(x_k),
    \qquad x_{k+1} = z_k - \beta_k^\Pi(z_k - \Pi_{g_{j_k}}z_k),
    \end{equation}
    with parameters ${\alpha_k> 0}$ and $ \beta_k^{\Pi} \in (0,2) $ for all ${k \in \N_0}$, and where $\Pi_{g_j}x$ denotes the projection onto $\{x\in \R^d:g_j(x) \leq 0\}$ for any $j \in \m$.
    \item Huber-penalty SGD (\citet[Eq. (4.1)]{nedich2023huber}), adapted to stochastic objective gradients
    \begin{equation}\label{alg:huber}
        x_{k+1} = x_k - s_k\bigl(
    \nabla f_{i_k}(x_k)
    + {\gamma_k^{\mathrm{H}}} \nabla h_{\delta_k^{\mathrm{H}}}(x_k;a_\jk,-b_\jk)
    \bigr),
    \end{equation}
    with positive parameters ${s_k}$, ${\gamma_k^{\mathrm{H}}}$, and ${\delta_k^{\mathrm{H}}}$ for all ${k \in \N_0}$.
    \item ABSGD \citep[Eq. (4)]{dimitrieski2025stochastic}
    \begin{equation}\label{alg:absgd}
        x_{k+1} = x_k - {\gamma_k^{\mathrm{AB}}}\Bigl(
    \nabla f_{i_k}(x_k)
    + \nabla B\bigl(g_{j_k}(x_k),\delta_k^{\mathrm{AB}}\bigr)
    \Bigr),
    \end{equation}
    with positive parameters ${\gamma_k^{\mathrm{AB}}}$ and ${ \delta_k^{\mathrm{AB}}}$ for all ${k \in \N_0}$.
    \item APDSG \citep[Algorithm 1 with Setting 2]{xu2020primal}
    \begin{equation}\label{alg:apdsg}
    \begin{aligned}[t]
        \psi_{j_k}(x_k)
        &= \max\{0,\varrho g_{j_k}(x_k) + \bar{z}_{k,j_k}\}
        \nabla g_{j_k}(x_k),\\
        x_{k+1}
        &= x_k - \mathbf{D}_k(\bar s_k)^{-1}\bigl(
        \nabla f_{i_k}(x_k) + \psi_{j_k}(x_k)
        \bigr),\\
        \bar{z}_{k+1,j}
        &= \begin{cases}
            \bar{z}_{k,j}, & j \neq j_k,\\
            \bar{z}_{k,j} + \bar{\rho}_k
            \max\left\{-\frac{\bar{z}_{k,j}}{\varrho},
            g_{j_k}(x_k)\right\}, & j = j_k,
        \end{cases}
    \end{aligned}
    \end{equation}
    where, for all ${k \in \N_0}$, we have $0 < \bar{\rho}_k \leq \varrho$ with $\varrho >0$, and ${\bar{z}_k \in \R^m}$ with~${\bar{z}_0 = 0}$. We define~${\mathbf{D}_k(\bar s_k) = \mathrm{diag}(\bar{s}_k) + \frac{1}{\bar{\alpha}_k}\mathrm{I}}$, with
    \(
    \bar{s}_k = \bar{\eta} \left(\sum_{t=0}^k \frac{( \nabla f_{i_t}(x_t) + \psi_{j_t}(x_t))^2}{\bar{\gamma}_t^2}\right)^{1/2},
    \) and $\bar{\alpha}_k > 0$ for all ${k \in \N_0}$,
    where~${\bar{\eta} > 0}$ and $\bar{\gamma}_t = \max\{1, \norm{ \nabla f_{i_t}(x_t) + \psi_{j_t}(x_t)}\}$ for every~${t \in \overline{0:k}}$.
\end{itemize}
\end{subequations}
In our setup, at $k \in \N_0$, all methods use i.i.d. uniformly sampled objective components (with indices $i_k \sim \Ucal(\n)$) and adaptively sampled constraints (with indices $j_k \in \m$), according to a probability distribution determined by the proposed adaptive CST in this paper, given in Algorithm~\ref{fig:CST} and Definition~\ref{def:CST:params}. The initial condition $x_0$ is treated as user-defined.

\section{Constraint-Sampling Technique}\label{sec:CST}
In this section we present an adaptive CST for solving Problem~\eqref{problem:1} with stochastic gradient-based methods such as~\eqref{alg:projection}--\eqref{alg:apdsg}.
The CST aims to enhance the performance of these methods by more frequently sampling recently violated constraints, and ultimately identify the active set~${\Acal\subseteq \m}$ in a time- and resource-efficient, online manner. In Section~\ref{sec:CST_explanation}, we first present the CST and then, in Section~\ref{sec:CST_theory}, provide further theoretical considerations.
\subsection{Algorithm  and Motivation}\label{sec:CST_explanation}

\begin{figure}[!t]
\centering
    \captionsetup{type=algorithm}

\definecolor{hfimp}{HTML}{0047AB}% draw from the important set (lower branch)
\definecolor{hfall}{HTML}{C2410C}% draw from all constraints (upper branch)
% Threshold of the test, set once for all four edge labels. The layout adapts to its length.
% For the drawn index, the vector is a_{j_k}:  \def\hfthr{\varepsilon_k\norm{a_{j_k}}}
\def\hfthr{-\varepsilon_k\norm{a_\jk}}
% measured in the font of the boxes: distance of two lines joined by \\[3pt], common width of
% the two boxes on the left, common width of the two draws, and the right column: the
% removal box with an inset \hfci on each side (at least 6.5mm, and enough for the label right
% of the arrow into its middle), the left inset holding the lane of the counter within its threshold
{\footnotesize
 \xdef\hfls{\the\dimexpr\baselineskip+3pt\relax}%
 \dimen0=0pt
 \def\hfmax#1{\settowidth{\dimen2}{#1}\ifdim\dimen2>\dimen0 \dimen0=\dimen2 \fi}%
 \hfmax{initialize}\hfmax{$x_0\in\R^d$}\hfmax{$k=0$, $\Gik=\emptyset$}%
 \hfmax{w.p.~$1-\eta$}\hfmax{sample}%
 \xdef\hflw{\the\dimexpr\dimen0+10pt\relax}%
 \settowidth{\dimen0}{$j_k\sim\mathcal{U}(\overline{1{:}m})$}%
 \settowidth{\dimen2}{$j_k\sim\mathcal{U}(\Gik)$}%
 \ifdim\dimen2>\dimen0 \dimen0=\dimen2 \fi
 \xdef\hfdw{\the\dimexpr\dimen0+10pt\relax}%
 \settowidth{\dimen0}{$\Gikk=\Gik\setminus\{j_k\}$}%
 \settowidth{\dimen2}{$u_{j_k,k+1}>t_\jk$}%
 \dimen2=\dimexpr\dimen2+1.4mm-(\dimen0+10pt)/2\relax
 \ifdim\dimen2<6.5mm \dimen2=6.5mm \fi
 \xdef\hfci{\the\dimen2}%
 \xdef\hfcw{\the\dimexpr\dimen0+10pt+2\dimen2\relax}}%
\scalebox{0.925}[0.925]{
\begin{tikzpicture}[
    font=\footnotesize,
    >={Stealth[length=2.1mm,width=1.7mm]},
    hfbox/.style={draw=black!65, rounded corners=1.5pt, line width=0.6pt,
        align=center, inner xsep=5pt, inner ysep=4pt, minimum height=7.5mm},
    hfimpbox/.style={hfbox, draw=hfimp, line width=0.9pt, fill=hfimp!20},
    hfallbox/.style={hfbox, draw=hfall, line width=0.9pt, fill=hfall!20},
    hfupd/.style={hfbox, fill=black!4},
    hfflow/.style={->, line width=0.7pt, draw=black!80},
    hfiflow/.style={hfflow, draw=hfimp, line width=0.9pt},
    hfaflow/.style={hfflow, draw=hfall, line width=0.9pt},
    hfguard/.style={inner sep=2pt},
]
% ---- insert/reset, and the draw from all constraints above it
\node[hfupd] (rs) {$\Gikk=\Gik\cup\{j_k\}$\\[2pt] $u_{j_k,k+1}=0$\\[2pt] $t_\jk=\tau_k$};
\node[hfallbox, minimum width=\hfdw, above=8mm of rs] (er) {$j_k\sim\mathcal{U}(\overline{1{:}m})$};
% ---- right column on one vertical axis, right of the label of the failed test: optimization
%      step in the row of the insert/reset, the removal below it, and the counter at the bottom
\node[hfguard, anchor=south west] (lErLt) at ($(er.east)+(2mm,0)$) {$g_{j_k}(x_k)<\hfthr$};
\coordinate (x0) at ($(lErLt.east)+(2mm,0)$);
\node[hfbox, anchor=west, minimum width=\hfcw] (opt) at (x0 |- rs) {optimization step with $g_{j_k}$};
\node[hfupd, anchor=north] (rm) at ($(opt.south)+(0,-5mm)$) {$\Gikk=\Gik\setminus\{j_k\}$};
\path let \p1=($(opt.east)-(opt.west)$) in
    node[hfupd, anchor=north, minimum width=\x1-0.6pt] (cnt) at ($(rm.south)+(0,-7mm)$)
    {$u_{j_k,k+1}=u_{j_k,k}+1$};
% ---- draw from the important set, in the row of the counter
\node[hfimpbox, minimum width=\hfdw] (ex) at (rs |- cnt) {$j_k\sim\mathcal{U}(\Gik)$};
\node[hfguard, anchor=north west] (lExLt) at ($(ex.east)+(2mm,0)$) {$g_{j_k}(x_k)<\hfthr$};
% ---- left column, both boxes of the same width: the sampling block (the condition on the
%      important set is on its input), whose lowest line is in the row of the draw from the
%      important set, and the start in the row of the draw from all constraints
\coordinate (sx) at ($(rs.west)+(-6.5mm,0)$);
\node[hfupd, minimum width=\hflw, anchor=mid east] (s1) at (sx |- ex)
    {\textcolor{hfall}{w.p.~$1-\eta$}\\[3pt] sample\\[3pt] \textcolor{hfimp}{w.p.~$\eta$}};
\node[hfbox, fill=white, line width=0.8pt, draw=black!80, minimum width=\hflw-0.2pt, anchor=east] (init)
    at (sx |- er)
    {initialize\\[2pt] $x_0\in\R^d$\\[2pt] $k=0$, $\Gikzero=\emptyset$};
% ---- edges: initialization and sampling. The start leads straight into the draw from all
%      constraints. The branch of the sampling block leaves at the height of its line and rises
%      along one line right of the left column into this line, and the branch of the empty
%      important set joins it below (with the loop back)
\draw[hfaflow] (init.east) -- (er.west);
\coordinate (xf) at ($(sx)+(3.25mm,0)$);
\coordinate (o1) at ($(s1.mid east)+(0,2*\hfls)$);
\draw[hfaflow, -, rounded corners=2.5pt] (o1) -- (xf |- o1) -- (xf |- er);
\draw[hfiflow] (s1.mid east) -- (ex.west);
% ---- edges: test
\draw[hfaflow] (er) -- node[right, hfguard, pos=0.5] (lUpGe) {$g_{j_k}(x_k)\ge\hfthr$} (rs);
\draw[hfiflow] (ex) -- node[right, hfguard, pos=0.7] (lLoGe) {$g_{j_k}(x_k)\ge\hfthr$} (rs);
\coordinate (xe) at ($(x0)+(8mm,0)$);
\draw[hfaflow, rounded corners=2.5pt] (er.east) -|
    node[left, hfguard, pos=0.75, text=black!65] (lNc) {} (opt.north -| xe);
\draw[hfiflow] (ex) -- (cnt);
% ---- insert/reset, counter within its threshold, and removal enter the optimization step
%      separately. The arrow into the removal box is vertical and enters it in the middle
\draw[hfflow] (rs) -- (opt);
\coordinate (xc) at ($(x0)+(0.5*\hfci,0)$);
\coordinate (ylab) at ($(cnt.north)!0.5!(rm.south)$);
\draw[hfiflow] (cnt.north -| xc) -- (opt.south -| xc);
\node[hfguard, anchor=east] (lCntLe) at ($(xc |- ylab)+(-0.3mm,0)$) {$u_{j_k,k+1}\le t_\jk$};
\draw[hfiflow] (cnt) -- (rm);
\node[hfguard, anchor=west] (lCntGt) at ($(rm.south |- ylab)+(0.3mm,0)$) {$u_{j_k,k+1}>t_\jk$};
\draw[hfiflow] (rm) -- (opt);
% ---- loop back just above the highest box or label. Under its condition, it enters the sampling
%      block from below when the important set is not empty. When it is empty, the loop back
%      joins the line into the draw from all constraints, midway between the two boxes on the left
\path let \p1=(init.north), \p2=(lErLt.north) in coordinate (top) at (0,{max(\y1,\y2)+3.5mm});
\coordinate (lft) at ($(init.west)+(-4.5mm,0)$);
\coordinate (ya) at ($(init.south)!0.5!(s1.north)$);
\draw[hfflow, rounded corners=3pt] (opt.north)
    -- node[right, hfguard, pos=0.3] (lLoop) {$k\gets k+1$} (opt.north |- top) -- (lft |- top)
    -- ($(lft |- s1.south)+(0,-6mm)$)
    -- node[above, hfguard, pos=0.5] (lGne) {$\Gik\neq\emptyset$} ($(s1.south)+(0,-6mm)$)
    -- (s1.south);
\draw[hfaflow, -] (lft |- ya) -- node[above, hfguard, pos=0.5] (lGe) {$\Gik=\emptyset$} (xf |- ya);
\end{tikzpicture}
}
\caption{Proposed adaptive constraint-sampling technique.}
\label{fig:CST}
\end{figure}

To start, we first define the setup of the proposed CST given in Algorithm~\ref{fig:CST} as follows.
\begin{definition}[Algorithm~\ref{fig:CST} Setup]
    \label{def:CST:params}
    Algorithm~\ref{fig:CST} is defined by 
    \begin{itemize}
        \item ${\eta \in (0,1)}$,
        \item $\{\tau_k\}_{k \in \N_0}$ with $\tau_0>0$, $\tau_{k+1} \geq\tau_{k}$,
        \item $\{\varepsilon_k\}_{k \in \N_0}$ with $\varepsilon_0 > 0$, $\varepsilon_k  > \varepsilon_{k+1}$, $\lim_{k \to \infty} \varepsilon_k = 0$,
    \end{itemize}
   and by
    \begin{itemize}
        \item ${\Gik \subseteq \overline{1:m}}$ for all $k \in \N_0$ with~${\G_{I,0} = \emptyset}$, 
        \item $(u_{j,k}, t_{j}) \in \mathbb N_0 \times \R_{> 0}$ for all $j \in \Gik$ and for $k \in \N_0$ where $\Gik \neq \emptyset$,
    \end{itemize}
    as determined in Algorithm~\ref{fig:CST}. 
\end{definition}
In particular, $\eta, \tau_k, \varepsilon_k$ are design parameters for the algorithm. They determine the adaptation of $\Gik \subseteq\m$, which is called the important set (of constraint indices) and is a maintained working set. The parameter ${\eta \in (0,1)}$ determines the prioritization of constraint indices in $\Gik$, and restricting  ${\eta < 1}$ ensures that all constraint indices are sampled with some positive probability. Additionally, $\{\tau_k\}_{k\in \N_0}$ is called the threshold sequence, and~$\{\varepsilon_k\}_{k\in \N_0}$ is called the 
(constraint) margin sequence. Finally,  for a constraint with index ${j \in \Gik \neq \emptyset}$, we call $u_{j,k} \in \N_0$ 
a counter that measures the corresponding constraint's 'lack of importance', and $t_j \in \R_{> 0}$ is a corresponding counter threshold. 
For example, ${u_{j,k}=\bar n}$ means that the constraint with index $j$ at iteration $k$ was consecutively (i.e., constantly) feasible with respect to the constraint margin in the last~${\bar n}$ draws from the sets $\mathbb{G}_{I,k_\ell}$, ${k_\ell < k}$, $\ell \in \overline{1:\bar n}$. If ${u_{j,k} > t_j}$,~then~the constraint with index~$j$ is considered
to be 'not important', and its index~$j$ is removed~from~$\Gik$.

The idea and motivation behind the design choices for Algorithm~\ref{fig:CST} is the following. Constraints that are violated should be sampled with priority, and therefore added to the important set. Constraints that are rarely violated are considered to be of less priority, and should not remain in the important set. To remove such a constraint with index $j$ from this set, consecutive 'non-violations' are counted
by~$u_{j,k}$. To guarantee that active constraints will eventually stay in the important set, violation is tested with respect to some margin~$\varepsilon_k$ and the threshold sequence $\{\tau_k\}_{k \in \N_0}$ is designed to be monotonically increasing.

\begin{figure}[t]
\centering
\definecolor{hfimp}{HTML}{0047AB}
\definecolor{hfall}{HTML}{C2410C}

\def\hfthr{-\varepsilon_k\|a_j\|}

\def\hftauoffset{0.625}
\pgfmathsetmacro{\hflevela}{0.75*sqrt(5)+\hftauoffset}
\pgfmathsetmacro{\hflevelb}{0.75*sqrt(15)+\hftauoffset}
\pgfmathsetmacro{\hflevelc}{0.75*sqrt(32)+\hftauoffset}
\begin{tikzpicture}[
    font=\footnotesize,
    hfmember/.style={fill=hfimp!20},
    hfcnt/.style={
        line width=1.2pt,
        draw=hfimp,
        line join=round
    },
    hfthr/.style={
        line width=1.1pt,
        draw=black,
        dash pattern=on 4pt off 2.2pt
    },
    hftau/.style={
        draw=black!60,
        line width=1pt,
        line cap=round,
        dash pattern=on 0pt off 3pt
    },
    hfguide/.style={
        draw=black!45,
        line width=0.45pt
    },
    hfx/.style={
        cross out,
        draw=black,
        line width=0.8pt,
        minimum size=5.5pt,
        inner sep=0pt
    },
    hfmV/.style={
        regular polygon,
        regular polygon sides=3,
        shape border rotate=180,
        fill=black,
        draw=black,
        inner sep=0pt,
        minimum size=8.5pt
    },
    hfmG/.style={
        circle,
        fill=hfimp,
        draw=hfimp,
        line width=0.6pt,
        inner sep=0pt,
        minimum size=6.6pt
    },
    hfmM/.style={
        circle,
        fill=hfall,
        draw=hfall,
        line width=0.6pt,
        inner sep=0pt,
        minimum size=6.6pt
    },
    hfaxis/.style={
        line width=0.65pt,
        draw=black!90
    },
    hftick/.style={
        font=\scriptsize,
        inner sep=1.5pt
    }
]
\begin{scope}[x=0.29cm, y=0.36cm]
  % Axes.
  \draw[hfaxis, -{Stealth[length=1.8mm,width=1.4mm]}]
    (0,0) -- (41.2,0) node[hftick, right] {$k$};
  \draw[hfaxis, -{Stealth[length=1.8mm,width=1.4mm]}]
    (0,0) -- (0,6.6);

  \foreach \k in {0,5,...,40}
    \draw[hfaxis] (\k,0) -- ++(0,-2.2pt)
      node[hftick, below] {$\k$};
  \foreach \v in {0,2,4,6}
    \draw[hfaxis] (0,\v) -- ++(-2.2pt,0)
      node[hftick, left] {$\v$};

  % Membership of j.
  \fill[hfmember] (5,0) rectangle (26,6.2);
  \fill[hfmember] (32,0) rectangle (40,6.2);

  % Guide lines for iterations with j_k = j.
  \foreach \k in {3,5,8,10,12,15,18,21,23,26,29,32,35,38}
    \draw[hfguide] (\k,6.8) -- (\k,0);

  % Grey dotted reference curve:
  % tau_k = (3/4)sqrt(k) + 3/4.
  \draw[hftau]
    plot[domain=0:40, samples=201, variable=\x]
      (\x,{0.75*sqrt(\x)+\hftauoffset});

  % Counter and threshold while j is important.
  \draw[hfcnt]
    (5,0) -- (8,0) -- (8,1) -- (12,1) -- (12,2)
    -- (15,2) -- (15,0) -- (18,0) -- (18,1)
    -- (21,1) -- (21,2) -- (23,2) -- (23,3) -- (26,3);

  \draw[hfthr]
    (5,\hflevela) -- (15,\hflevela)
    -- (15,\hflevelb) -- (26,\hflevelb);

  \draw[hfcnt]
    (32,0) -- (35,0) -- (35,1)
    -- (38,1) -- (38,2) -- (40,2);

  \draw[hfthr]
    (32,\hflevelc) -- (40,\hflevelc);

  % Keep the original counter extension and cross.
  \draw[hfcnt] (26,3) -- (26,4);
  \node[hfx] (rm) at (26,4) {};

  % Iterations with j_k = j; the resulting change appears at k+1.
  \foreach \k in {5,15,32}
    \node[hfmV] at (\k,7.15) {};
  \foreach \k in {8,12,18,21,23,26,35,38}
    \node[hfmG] at (\k,7.15) {};
  \foreach \k in {3,10,29}
    \node[hfmM] at (\k,7.15) {};
\end{scope}

% Legend (two rows).
\matrix (leg1) [
    anchor=north,
    inner sep=0pt,
    column sep=1.3mm,
    nodes={
        font=\scriptsize,
        anchor=west,
        inner sep=1pt
    }
] at (5.95,-0.72) {
  \node[hfmV] {}; &
  \node {$g_j(x_k)\ge\hfthr$}; &
  \node {}; &
  \node[hfmG] {}; &
  \node {$j_k\sim\mathcal{U}(\Gik)$,
    \ $g_j(x_k)<\hfthr$}; &
  \node {}; &
  \node[hfmM] {}; &
  \node {$j_k\sim\mathcal{U}(\overline{1{:}m})$,
    \ $g_j(x_k)<\hfthr$}; \\
};

\matrix [
    anchor=north,
    inner sep=0pt,
    column sep=1.3mm,
    nodes={
        font=\scriptsize,
        anchor=west,
        inner sep=1pt
    }
] at ($(leg1.south)+(0,-1.6mm)$) {
  \draw[hfcnt] (0,0) -- (5mm,0); &
  \node {counter $u_{j,k}$}; &
  \node {}; &
  \draw[hfthr] (0,0) -- (5mm,0); &
  \node {threshold $t_j$}; &
  \node {}; &
  \draw[hftau] (0,0) -- (5mm,0); &
  \node {$\tau_k$}; &
  \node {}; &
  \fill[hfmember] (0,-1.1mm) rectangle (5mm,1.1mm); &
  \node {$j\in\Gik$}; \\
};
\end{tikzpicture}
\caption{Illustrative example for Algorithm~\ref{fig:CST}: Sample path for one fixed
constraint with index~${j\in\m}$ plotted over the iterations $k\in \N_0$. Depicted: membership in $\Gik$ (shaded blue means~${j\in\Gik}$),
counter $u_{j,k}$ (blue line), threshold $t_j$ (dashed line),
and $\tau_k$ (dotted line), with
$\tau_k=\tfrac{3}{4}\sqrt{k}+\tfrac{5}{8}$ for $k\in \N_0$,
chosen for convenience.
A vertical line at an instant~$k$ indicates that the $j$-th constraint
is sampled (either from $\Gik$  or $\overline{1:m}$), i.e., $j_k=j$.
A black triangle at the top of the vertical line 
marks that the constraint is violated and results in
an insertion into $\Gikk$ or a reset of its counter with an update of its threshold,
i.e., $j\in\Gikk$, $u_{j,k+1}=0$ and~${t_j=\tau_k}$.
A blue circle marks that the constraint is satisfied and
an increment $u_{j,k+1}=u_{j,k}+1$, after which $j$ is removed if and only if
$u_{j,k+1}>t_j$.
A red circle marks that the constraint is satisfied  and
an iteration in which none of the three
quantities changes, even if $j\in\Gik$, as, e.g., at $k=10$.
Further, for instance, after the reset (black triangle)
at $k=15$, which sets $t_j=\tau_{15}\approx3.530$,
the increments at $k=18$, $21$, $23$, and~$26$ raise the counter to
$u_{j,27}=4>t_j$ (cross).
Thus, $j$ is removed at~$k=26$, i.e., ${j\in\mathbb{G}_{I,26}}$ and $j\notin\mathbb{G}_{I,27}$.}
\label{fig:sample-path}
\end{figure}
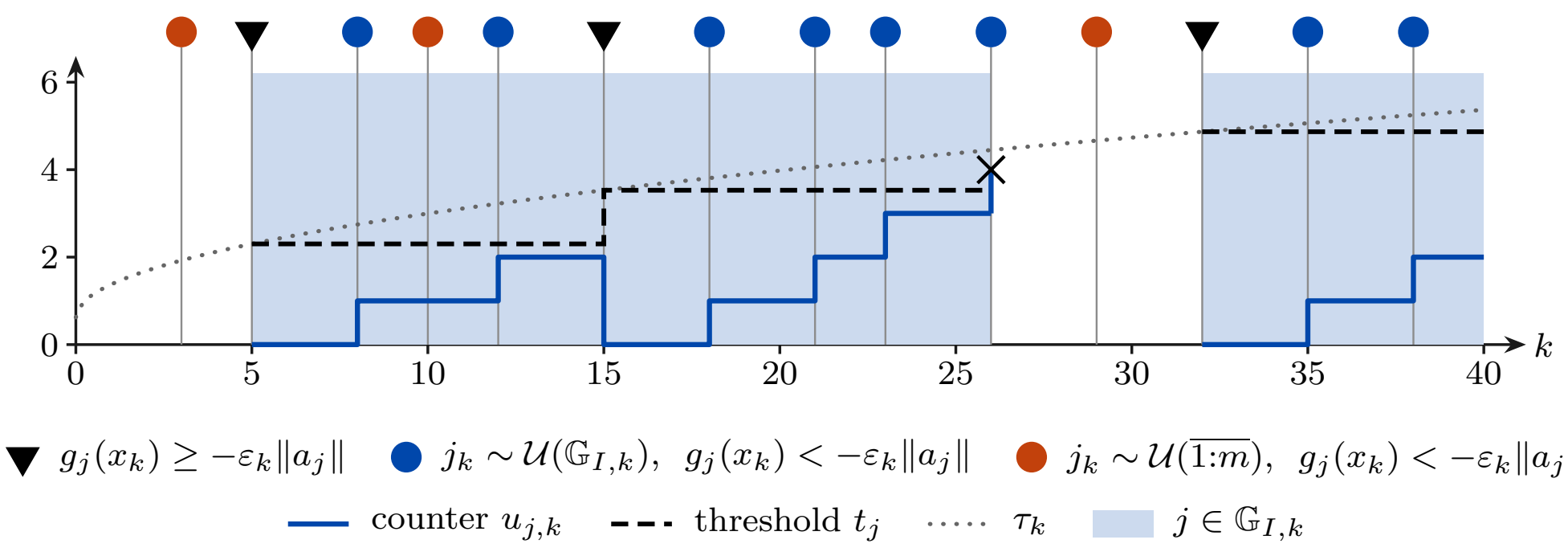

In the following, we describe in detail the flowchart of Algorithm~\ref{fig:CST}.
At iteration ${k=0}$, the important set is initialized with ${\mathbb{G}_{I,0}=\emptyset}$ and a constraint with index $j_k$ is sampled uniformly from $\m$. If this constraint passes the margin-based 'violation' test, $j_k$ is added to the important set, its counter and threshold are initialized, and an optimization step with the sampled constraint is performed. Otherwise, only the optimization step with the sampled constraint is performed. At any following iteration $k+1$, we first check if~${\mathbb{G}_{I,k+1} = \emptyset}$. If this holds, then the previous procedure is repeated. Further, if ${\mathbb{G}_{I,k+1} \neq \emptyset}$, with probability $1- \eta$, the previous procedure is repeated as well. If ${\mathbb{G}_{I,k+1} \neq \emptyset}$, with probability $ \eta$,
a constraint with index $j_{k+1}$ is sampled uniformly from the important set. If this constraint passes the 'violation' test, then $j_{k+1}$ is kept in the important set $\mathbb{G}_{I,k+2}$, its counter $u_{j_{k+1},k+2}$ is reset to zero and its threshold $t_{j_{k+1}}$ is updated to $\tau_{k+1}$, and an optimization step with the sampled constraint is performed. Otherwise, its counter is incremented. Subsequently, it is checked if ${u_{j_{k+1},k+2} > t_{j_{k+1}}}$, and in case this holds,~$j_{k+1}$ is removed from the important set $\mathbb{G}_{I,k+2}$. Thereafter, an optimization step with the sampled constraint is performed, and we proceed to the next iteration. For illustration, an example sample path for a fixed index~${j \in \m}$ is given in Figure~\ref{fig:sample-path}.

Notice that Algorithm~\ref{fig:CST} can be implemented with a per-iteration arithmetic overhead of either~$\mathcal{O}(1)$, or of amortized~$\mathcal{O}(1)$, depending on user preference and available memory resources. However, an implementation-specific discussion is out of the scope of this~paper. 

\subsection{Theoretical Considerations}\label{sec:CST_theory}
While a rigorous convergence analysis is out of the scope of this paper due to space limitations and is deferred to future work, we provide informal reasoning as to why, for suitably chosen design parameters, Algorithm~\ref{fig:CST} identifies the active set for the methods~\eqref{alg:projection}--\eqref{alg:absgd}.
For this, we make the following standard assumption on the objective and the constraints of Problem~\eqref{problem:1} (see, e.g., \citet{boyd2004convex,nesterov2018convex}).
\begin{assumption}\label{ass:objective}
    Consider Problem \eqref{problem:1}. For every~$f_i \in \mathcal{C}^1$,~$i \in \n$, and for every~${x,y \in \R^d}$, we assume that
    \begin{enumerate}[label=\textup{\alph*)}]
        \item \label{ass:lsmoothe}  $f_i$ is $L$-smooth, i.e., $\|\nabla f_i(y) - \nabla f_i(x)\| \leq L\|y-x\|$,
        \item \label{ass:coerce} $f$ is $\mu$-strongly convex, i.e., $f(y) - f(x) \geq \nabla f^\top(x)(y-x) + \tfrac{\mu}{2} \|y-x\|^2$.
    \end{enumerate}
    Moreover, we assume a non-empty interior of ${\mathbb{X} \coloneqq \left\{ x \in \mathbb{R}^d \colon g_j(x) \leq 0, \ j \in \overline{1:m} \right\}}$.
\end{assumption}

Notice that, by the construction of Algorithm~\ref{fig:CST}, every constraint with index $j \in \m$ is sampled infinitely often. Then, for methods \eqref{alg:projection}--\eqref{alg:absgd} under Assumption~\ref{ass:objective}, and with a suitable hyperparameter choice, it can be shown that the iterates almost surely converge to~$x^*$ with some rate. If $\{\varepsilon_k\}_{k\in \N_0}$ is then chosen to asymptotically vanish at a slower rate than~$\|x_k-x^{*}\|$, the active constraints will eventually be retained in the important set, and all non-active ones eventually removed from it. This can be achieved since each constraint in the important set is, by design of Algorithm~\ref{fig:CST}, sampled with some positive probability.

\section{Numerical Results}\label{sec:empirical_results}

In this section, we present numerical results for Algorithm~\ref{fig:CST}, implemented for different stochastic gradient-based methods. We first demonstrate the effect of Algorithm~\ref{fig:CST} and its scalability on synthetic optimization problems. Then, we apply it to the problem of learning a CBF, as adapted from \cite{lutkus2025incremental}.

\subsection{Synthetic Optimization Problems}\label{section:numerical_results}

\textbf{Setup:} We consider a~${d = 100}$-dimensional problem, where the objective is a sum of~\({n_Q = 5}\) quadratic and \(n_{\log} = 5 \) logarithmic components of form
\begin{equation*}
        \label{problem:section:numerical_results}
        f(x) \coloneqq \tfrac{1}{n_Q + n_{\log}}\sum\nolimits_{l=1}^{n_Q} (x-\nu_l)^\top Q_l(x-\nu_l) + \tfrac{C_{\log}}{n_Q + n_{\log}}\sum\nolimits_{\ell=1}^{n_{\log}} \log\Bigl(1+\exp\bigl(\kappa_\ell^\top(x-\chi_\ell)\bigr)\Bigr),
\end{equation*}
where \(Q_l \in \R^{d\times d}\), \(Q_l \succ 0\), \( \nu_l \in \R^d\), \(l\in \overline{1:n_Q}\), and \( C_{\log} > 0 \), \(\kappa_\ell \in \R^d\), \(\chi_\ell \in \R^d\), \({\ell \in \overline{1:n_{\log}}}\).
The constraints are chosen 
to be affine, i.e.,~${g_j(x) \coloneqq a_j^\top x + b_j}$, with \(a_j \in \mathbb{R}^d \) and~\({b_j \in \mathbb{R}}\)
for~${j\in \m}$ and generated in a random manner. We consider problems with three comparatively quite different constraint configurations, and same constrained minimizer~$x^*$:
\begin{enumerate}[
    label={4.1~\roman*)},
    ref={4.1~\roman*)},
    align=left,
    labelindent=10pt,
    leftmargin=*,
    widest=iii
]
    \item \label{problem:type:1}
    a baseline configuration with $2 \cdot 10^5$ constraints,

    \item \label{problem:type:2}
    a large-scale configuration with $10^6$ constraints,

    \item \label{problem:type:3}
    Problem~\ref{problem:type:1} with $10^3$ constraints added to the immediate proximity
    of $x^*$.
\end{enumerate}
We apply Algorithm~\ref{fig:CST} with the methods \eqref{alg:projection}--\eqref{alg:apdsg}. We implement it both with a single sampled constraint and a constraint mini-batch (CMB) of size $8$, and benchmark it against uniform constraint sampling.
For the CMB approach, we sequentially sample a constraint with replacement and update the important set, and use the averaged gradient information with fixed primal (and dual) variables for the update step. For Projection-SGD \eqref{alg:projection}, we use the averaged projected steps as in \citet[Algorithm~1]{wang2015random}.
Each method's hyperparameters and the hyperparameters of Algorithm~\ref{fig:CST} are given in the appendix.

\begin{figure*}[t]
    \centering
\scalebox{1}{
    % ---------- Row 1 --------
    \subfigure[Problem~\ref{problem:type:1}]{\label{fig:200_0_1}%
        \includegraphics[width=0.32\textwidth]{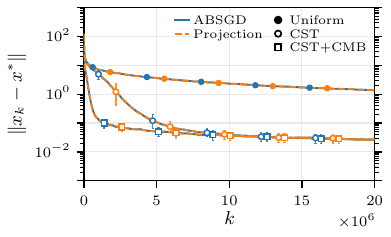}}%
    \hfill
    \subfigure[Problem~\ref{problem:type:2}]
    % [ABSGD \& Projection-SGD for problem~\ref{problem:type:2}]
    {\label{fig:1000_0_1}%
        \includegraphics[width=0.32\textwidth]{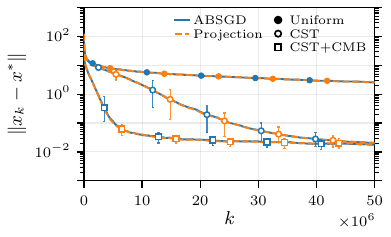}}%
    \hfill
    \subfigure[Problem~\ref{problem:type:3}]
    % [ABSGD \& Projection-SGD for problem~\ref{problem:type:3}]
    {\label{fig:201_0_1}%
        \includegraphics[width=0.32\textwidth]{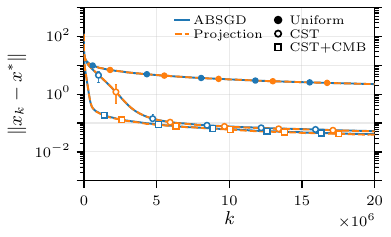}}
}
\scalebox{1}{
    % ---------- Row 2 ----------
    \subfigure[Problem~\ref{problem:type:1}]
    % [Huber-SGD \& APDSG for problem~\ref{problem:type:1}]
    {\label{fig:200_0_2}%
        \includegraphics[width=0.32\textwidth]{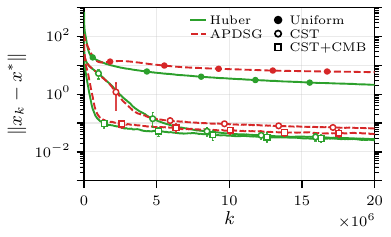}}%
    \hfill
    \subfigure[Problem~\ref{problem:type:2}]
    % [Huber-SGD \& APDSG for problem~\ref{problem:type:2}]
    {\label{fig:1000_0_2}%
        \includegraphics[width=0.32\textwidth]{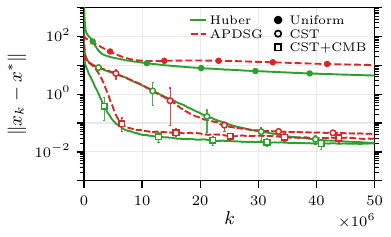}}%
    \hfill
    \subfigure[Problem~\ref{problem:type:3}]
    % [Huber-SGD \& APDSG for problem~\ref{problem:type:3}]
    {\label{fig:201_0_2}%
        \includegraphics[width=0.32\textwidth]{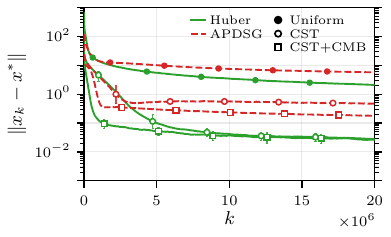}}
}
    \caption{Simulation results for Problems \ref{problem:type:1}--\ref{problem:type:3} with methods \eqref{alg:projection}--\eqref{alg:apdsg} and Algorithm~\ref{fig:CST}, implemented according to Definition~\ref{def:CST:params}. Mean Euclidean distance of the iterates to the constrained minimizer is plotted over the iteration count~${k \in \N_0}$ from $101$ sample trajectories. Error bars show $10 \%$ and $90 \%$ quantiles.
    }
    \label{fig:x_conv}
\end{figure*}
\begin{figure*}[t]
    \centering
    \scalebox{1}{
    \subfigure[\mbox{${\eta \in \{0.1, 0.5,0.9\}}$}, \mbox{${\omega = 0.6}$}]{\label{fig:set_200_1}%
        \includegraphics[width=0.32\textwidth]{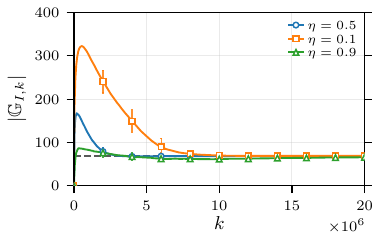}}%
    \hfill
    \subfigure[\mbox{${\eta = 0.5}$}, \mbox{${\omega \in \{0.5, 0.6, 0.9\}}$}]{\label{fig:set_200_2}%
        \includegraphics[width=0.32\textwidth]{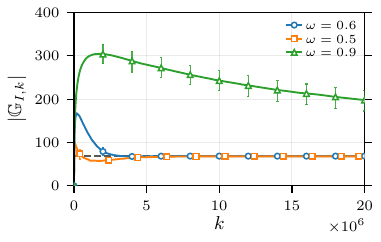}}%
    \hfill
    \subfigure[\mbox{$(\eta,\omega) = (0.5,0.6)$}, \mbox{${\mathrm{CMB}\in \{1,4,8\}}$}]{\label{fig:set_200_3}%
        \includegraphics[width=0.32\textwidth]{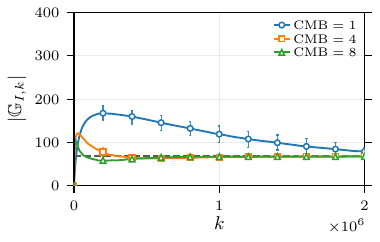}}%
    \hfill
    }
    \caption{Simulation results for Problem \ref{problem:type:1} with ABSGD \eqref{alg:absgd} and Algorithm~\ref{fig:CST}, implemented according to Definition~\ref{def:CST:params}.
    Mean cardinality of the set~$\Gik$ is plotted over the iteration count~${k \in \N_0}$ from $101$ sample trajectories.
    Error bars show $10 \%$ and~$90 \%$ quantiles. Dashed line denotes the number of active constraints.
    }
    \label{fig:sets}
\end{figure*}
\noindent\textbf{Results:} Firstly, we analyze the iterate convergence in Figure~\ref{fig:x_conv}, plotted over the iteration count $k \in \N_0$. It can be observed that across all methods and problem instances, Algorithm~\ref{fig:CST} decreases the iterates' final measured Euclidean distance to~$x^*$ on average by two orders of magnitude for ABSGD (\eqref{alg:absgd}, blue), Projection-SGD (\eqref{alg:projection}, orange), and Huber-penalty SGD (\eqref{alg:huber}, green), and by one to two orders of magnitude for APDSG (\eqref{alg:apdsg}, red). Adding a CMB improves the initial performance. Comparing the results for Problem~\ref{problem:type:2} (Figure~\ref{fig:1000_0_1} and~\ref{fig:1000_0_2}) to the results for Problem~\ref{problem:type:1} (Figure~\ref{fig:200_0_1} and~\ref{fig:200_0_2}) shows that Algorithm~\ref{fig:CST} remains effective when increasing~$m$. The results for Problem~\ref{problem:type:3} (Figure~\ref{fig:201_0_1} and~\ref{fig:201_0_2}) are similar to the results for Problem~\ref{problem:type:1}. Overall, compared to uniform constraint sampling, Algorithm~\ref{fig:CST} gives a significant performance improvement, which can be over $100$ times better.

Secondly, we analyze the dynamics of Algorithm~\ref{fig:CST} for different hyperparameter configurations in Figure~\ref{fig:sets}. For this, we plot the cardinality of~$\Gik$ for ABSGD \eqref{alg:absgd} on Problem~\ref{problem:type:1} over the iteration count $k \in \N_0$. For the baseline configuration with~${\eta = 0.5}$ and~${\tau_k = (k+1)^{\omega}}$ with $\omega = 0.6$, ${k \in \N_0}$, (Figure~\ref{fig:set_200_1}--\ref{fig:set_200_3}, blue), we observe that~$|\Gik| \gg |\Acal|$ initially, and then~$|\Gik|$ reduces to $|\Acal|$. We have verified that then also $\Gik = \Acal$ holds. Increasing~$\eta\in (0,1)$ results in a smaller initial overshoot (Figure~\ref{fig:set_200_1}, green). However, a subsequent undershoot is observed. Similar effects are observed for decreasing $\omega>0$ (Figure~\ref{fig:set_200_2}, orange). Moreover, increasing the CMB size leads to significantly faster convergence of $|\Gik|$ to $|\Acal|$ (Figure~\ref{fig:set_200_3}).

\noindent \textbf{Remark:} By comparing the wall-clock times of uniform constraint sampling and constraint sampling according to Algorithm~\ref{fig:CST}, we observed an increase of less than 10$\%$ across all simulations, regardless of the number of constraints.

\subsection{Learning a Control Barrier Function}\label{sec:dubins}
We aim to learn a CBF for a system of the form $\dot{x} = F(x) + G(x)u$, 
where $F:\R^{d_x} \to \R^{d_x}$ and $G:\R^{d_x} \to \R^{d_x \times d_u}$ are locally Lipschitz continuous,~${x \in \R^{d_x}}$ is the state, and~${u \in \R^{d_u}}$ is the control input. In particular, we consider the Dubins Car setup from \citet{lutkus2025incremental}. Following \citet{lutkus2025incremental}, we define a CBF for such a system as follows, where $\text{int}(S)$ denotes the interior and $\text{bd}(S)$ the boundary of some set~$S$:
\begin{definition}[\citet{lutkus2025incremental}]
Given a compact set $S \subseteq \R^{d_x}$, a set of control inputs $U \subseteq \R^{d_u}$, and a locally Lipschitz continuous extended class-$\Kcal$ function $\varphi:\R\to\R$, a continuously differentiable function $\bar{h}:D \to \R$ is said to be a control barrier function on an open set $D \supseteq S$ if:
    \begin{align*}
    &\bar{h}(x) > 0 \ \forall x \in \mathrm{int}(S), \qquad \bar{h}(x) < 0\ \forall x \in D \backslash S, \qquad \bar{h}(x) = 0 \ \forall x \in \mathrm{bd}(S)\\
    &\sup_{u \in U} \nabla \bar{h}^\top(x) \bigl(F(x)+G(x)u\bigr) \geq -\varphi\bigl(\bar{h}(x)\bigr) \ \forall x \in D.
    \end{align*}
    
\end{definition}

To learn a candidate for such a CBF from data, one approach is to parameterize the CBF as~\(
    \bar{h}_\theta(x)\coloneqq \theta^\top \phi(x) - b
\)
with $N$ basis functions defined by $\phi :\R^{d_x} \to \R^N$, with parameter~${\theta \in \R^N}$ and some ${b > 0}$, and solve the following constrained optimization problem to find a suitable~${\theta \in \R^N}$ (see, e.g., \citet{lutkus2025incremental})
\begin{equation}\label{problem:dubins}
\begin{aligned}
    &\underset{\theta \in \R^N}{\mathrm{min}} \hspace{0.1cm} \|\theta\|^2 \hspace{0.4cm}\\
    &\underset{\phantom{\theta \in \R^N}}{\mathrm{s.t.}} \ \bar{h}_\theta(x_i) \geq \epsilon_{\text{safe}} \hspace{0.1cm} \forall (x_i,u_i) \in \Scal, \qquad \qquad \bar{h}_\theta(x_i) \leq -\epsilon_{\text{unsafe}} \hspace{0.1cm}\forall (x_i,u_i) \in \bar{\Dcal},\\
    & \phantom{\underset{\phantom{\theta \in \R^N}}{\mathrm{s.t.}}  } \ \nabla \bar{h}_\theta^\top(x_i)\bigr(F(x_i) + G(x_i)u_i\bigl) + \varphi\bigl(\bar{h}_\theta(x_i)\bigr) \geq \epsilon_{\text{dyn}} \ \forall (x_i,u_i) \in \Dcal_{\mathrm{aug}},
    \end{aligned}
\end{equation}
where $\epsilon_{\text{safe}},\epsilon_{\text{unsafe}},\epsilon_{\text{dyn}} >0$ and $b > 0$ are hyperparameters. Additionally, ${\Dcal\coloneqq \{(x_i,u_i)\}_{i=1}^T}$ is a dataset of $T \in \mathbb N$ state-action pairs of observed system behavior and ${\Scal \subseteq \Dcal}$ is a so-called safe dataset (with some margin). The dataset~$\bar{\Dcal}$ consists of state-action pairs considered unsafe, and~$\Dcal_{\mathrm{aug}}$ is an augmented dataset based on $\Dcal$, designed to have sufficiently many state-action pairs for the training process.

\noindent\textbf{Setup:} We consider local optimization problems of form~\eqref{problem:dubins} adapted from \citet{lutkus2025incremental}, denoted by Dubins Car 4.2~i)--iv), and apply Algorithm~\ref{fig:CST} with the methods~\eqref{alg:projection} and~\eqref{alg:absgd}. We implement it both with a single sampled constraint and a CMB of size~$8$, and benchmark it against uniform constraint sampling. We use compactly supported radial basis functions following \citet{lutkus2025incremental}, slightly modified by replacing the angular-distance term with a chordal distance (see, e.g., \citet{flyer2009aradial}), to ensure periodicity and continuity of $\bar h_\theta$ and $\nabla \bar h_\theta$ in the angular terms. We set $\varphi$ to the identity function.
Additionally, we augment the objective by adding the term~${\sum_{(x_i,u_i) \in \bar{\Dcal} \cup \Dcal_\mathrm{aug}} \|\nabla \bar{h}_\theta(x_i)\|^2}$, as is done in several implementations (see, e.g., \citet{robey2020learning, lutkus2025incremental}), and normalize it by $|\bar{\Dcal} \cup \Dcal_\mathrm{aug}|$, and $\|\theta\|^2$ by $N$.
Each method's hyperparameters, given in the appendix together with the problem instance hyperparameters, are tuned to achieve the smallest Euclidean distance to a reference solution~$\theta^* \in \R^N$.
\begin{figure*}[t]
    \centering
    % ---------- Row 1 --------
    \subfigure[Dubins Car 4.2~i)]{\label{fig:dubins0}%
        \includegraphics[width=0.4\textwidth]{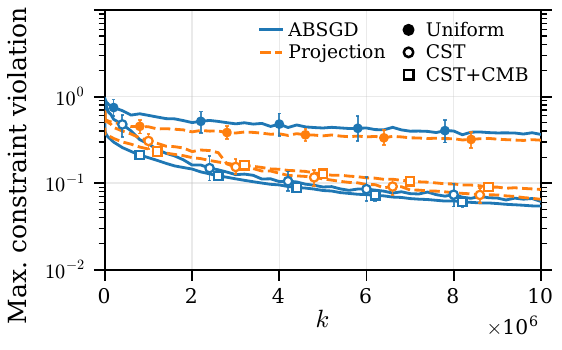}}%
    \hfill
    \subfigure[Dubins Car 4.2~ii)]{\label{fig:dubins1}%
        \includegraphics[width=0.4\textwidth]{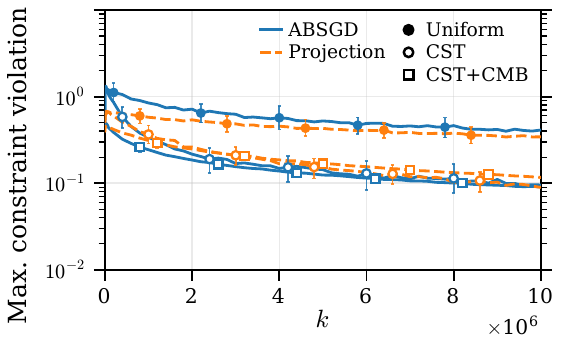}}%
   
    % ---------- Row 2 ----------
    % ---------- Row 1 --------
    \subfigure[Dubins Car 4.2~iii)]{\label{fig:dubins2}%
        \includegraphics[width=0.4\textwidth]{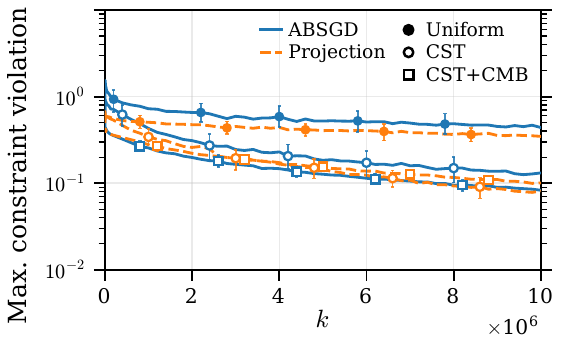}}%
    \hfill
    \subfigure[Dubins Car 4.2~iv)]{\label{fig:dubins3}%
        \includegraphics[width=0.4\textwidth]{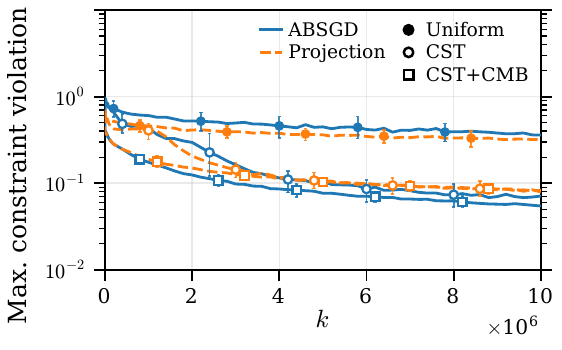}}%
        
    \vspace{-0.2cm}
    
    \caption{Simulation results for Dubins Car 4.2 i)--iv) with methods~\eqref{alg:projection} and~\eqref{alg:absgd}, implemented with Algorithm~\ref{fig:CST}, with setup as in Definition~\ref{def:CST:params}. Performance metric is the magnitude of the most violated constraint of Problem~\eqref{problem:dubins} with augmented objective, evaluated at every iteration~${k\in \N_0}$. Plotted: Mean maximum constraint violation over the iteration count~${k \in \N_0}$ from~$100$ sample trajectories. Error bars show~$10 \%$ and $90 \%$ quantiles.
    }
    \label{fig:dubins}
\end{figure*}

\noindent\textbf{Results:} We find a numerical solution of Problem~\eqref{problem:dubins} (with augmented objective) to obtain a CBF candidate, whose validity as a CBF is primarily determined by the violation of the constraints of this problem. Hence, for evaluation, we analyze the maximum constraint violation of the CBF candidate, as shown in \Cref{fig:dubins}, plotted over the iteration count $k \in \N_0$. In particular, across all problems and both methods, Algorithm~\ref{fig:CST} reduces the maximum constraint violation on average by half an order and up to one order of magnitude compared to uniform constraint sampling at the final iteration. Therefore, the CBF candidate obtained with Algorithm~\ref{fig:CST} is comparatively a significantly better candidate than the one obtained with uniform constraint sampling after $10^7$ iterations, with up to $10$ times smaller maximum constraint violation. Including a CMB further improves the initial performance, as depicted in \Cref{fig:dubins}.

\section{Conclusions and Outlook}
In this paper, we presented an adaptive CST for stochastic gradient-based methods that adds relatively minor per-iteration arithmetic overhead which is independent of the number of constraints. Numerical results show a potential performance improvement of two orders of magnitude and even more over uniform constraint sampling, especially in large-scale problems.

Theoretical results for the proposed constraint-sampling technique will be presented in an upcoming paper, in which both, almost sure asymptotic iterate convergence and finite-time active-set identification will be analyzed. For future work, it is interesting to, for example, investigate the performance of Algorithm~\ref{fig:CST} on non-linearly constrained problems, such as a (Deep) Neural Network approach for the CBF \citep{robey2020learning}, or with momentum-based methods or gradient-free methods (e.g., \citet{dimitrieski2025time}).

\bibliography{references}

@inproceedings{lutkus2025incremental,
  author={Lutkus, Paul and Anantharaman, Deepika and Tu, Stephen and Lindemann, Lars},
  booktitle={64{th} Conference on Decision and Control}, 
  title={Incremental composition of learned control barrier functions in unknown environments}, 
  year={2025},
  pages={5967-5972},
publisher = {IEEE}
}

@article{wang2016stochastic,
author = {Wang, Mengdi and Bertsekas, Dimitri P.},
title = {Stochastic first-order methods with random constraint projection},
journal = {SIAM Journal on Optimization},
volume = {26},
pages = {681--717},
year = {2016},
publisher={SIAM}
}

@article{dimitrieski2025stochastic,
  author={Dimitrieski, Naum and Cao, Jing and Ebenbauer, Christian},
  journal={IEEE Control Systems Letters}, 
  title={Stochastic gradient descent for constrained optimization based on adaptive relaxed barrier functions}, 
  volume={9},
  pages={829--834},
  year={2025},
  publisher={IEEE}
  }

@article{nedich2023huber,
      title={Huber loss-based penalty approach to problems with linear constraints}, 
      author={Angelia Nedi{\'c} and Tatiana Tatarenko},
        journal = {arXiv preprint arXiv:2311.00874},
      year={2023},
}

@article{xu2020primal,
author = {Xu, Yangyang},
title = {Primal-dual stochastic gradient method for convex programs with many functional constraints},
journal = {SIAM Journal on Optimization},
volume = {30},
pages = {1664--1692},
year = {2020},
publisher={SIAM}
}

@article{gallego2025cooper,
  title={Cooper: A library for constrained optimization in deep learning},
  author={Gallego-Posada, Jose and Ramirez, Juan and Hashemizadeh, Meraj and Lacoste-Julien, Simon},
  journal={arXiv preprint arXiv:2504.01212},
  year={2025}
}

@article{robey2021adversarial,
  title={Adversarial robustness with semi-infinite constrained learning},
  author={Robey, Alexander and Chamon, Luiz and Pappas, George J. and Hassani, Hamed and Ribeiro, Alejandro},
  journal={Advances in Neural Information Processing Systems},
  volume={34},
  pages={6198--6215},
  year={2021}
}

@article{rezkellah2025machine,
  title={Machine unlearning meets adversarial robustness via constrained interventions on {LLM}s},
  author={Rezkellah, Fatmazohra and Dakhmouche, Ramzi},
  journal={arXiv preprint arXiv:2510.03567},
  year={2025}
}

@article{de2003linear,
  title={The linear programming approach to approximate dynamic programming},
  author={De Farias, Daniela Pucci and Van Roy, Benjamin},
  journal={Operations Research},
  volume={51},
  pages={850--865},
  year={2003},
  publisher={INFORMS}
}

@article{de2004constraint,
  title={On constraint sampling in the linear programming approach to approximate dynamic programming},
  author={De Farias, Daniela Pucci and Van Roy, Benjamin},
  journal={Mathematics of Operations Research},
  volume={29},
  pages={462--478},
  year={2004},
  publisher={INFORMS}
}

@article{samiedaluie2024approximate,
  title={Approximate linear programming for a queueing control problem},
  author={Samiedaluie, Saied and Zhang, Dan and Zhang, Rui},
  journal={Computers \& Operations Research},
  volume={169},
  pages={106711},
  year={2024},
  publisher={Elsevier}
}

@inproceedings{cotter2016light,
  title={A light touch for heavily constrained {SGD}},
  author={Cotter, Andrew and Gupta, Maya and Pfeifer, Jan},
  booktitle={Conference on Learning Theory},
  pages={729--771},
  year={2016},
  organization={PMLR}
}

@book{boyd2004convex,
  title     = {Convex optimization},
  author    = {Boyd, Stephen and Vandenberghe, Lieven},
  year      = {2004},
  publisher = {Cambridge University Press}
}

@inproceedings{dimitrieski2025time,
  author={Dimitrieski, Naum and Reyer, Michael and Belabbas, Mohamed-Ali and Ebenbauer, Christian},
  booktitle={2025 European Control Conference}, 
  title={Time-delay induced stochastic optimization and extremum seeking}, 
  year={2025},
  pages={1636--1641},
publisher = {IEEE}
}

@inproceedings{robey2020learning,
  author={Robey, Alexander and Hu, Haimin and Lindemann, Lars and Zhang, Hanwen and Dimarogonas, Dimos V. and Tu, Stephen and Matni, Nikolai},
  booktitle={59{th} Conference on Decision and Control}, 
  title={Learning control barrier functions from expert demonstrations}, 
  year={2020},
  pages={3717--3724},
  publisher = {IEEE}
}

@book{nesterov2018convex,
  title={Lectures on convex optimization},
  author={Nesterov, Yurii},
  volume={137},
  series={Springer Optimization and Its Applications},
  year={2018},
  edition = {2nd},
  publisher={Springer Science \& Business Media}
}

@article{calafiore2006scenario,
  title={The scenario approach to robust control design},
  author={Calafiore, Giuseppe Carlo and Campi, Marco C.},
  journal={IEEE Transactions on Automatic Control},
  volume={51},
  pages={742--753},
  year={2006},
  publisher={IEEE}
}

@article{campi2009scenario,
  title={The scenario approach for systems and control design},
  author={Campi, Marco C. and Garatti, Simone and Prandini, Maria},
  journal={Annual Reviews in Control},
  volume={33},
  pages={149--157},
  year={2009},
  publisher={Elsevier}
}

@article{flyer2009aradial,
    author = {Flyer, Natasha and Wright, Grady B.},
    title = {A radial basis function method for the shallow water equations on a sphere},
    journal = {Proceedings of the Royal Society A: Mathematical, Physical and Engineering Sciences},
    volume = {465},
    pages = {1949--1976},
    year = {2009},
    publisher = {Royal Society Publishing}
}

@inproceedings{kotary2021ijcai,
  title     = {End-to-end constrained optimization learning: A survey},
  author    = {Kotary, James and Fioretto, Ferdinando and Van Hentenryck, Pascal and Wilder, Bryan},
  booktitle = {30{th} International Joint Conference on
               Artificial Intelligence},
  publisher = {International Joint Conferences on Artificial Intelligence Organization},
  pages     = {4475--4482},
  year      = {2021}
}

@article{laiu2019aconstraint,
  title={A constraint-reduced {MPC} algorithm for convex quadratic programming, with a modified active set identification scheme},
  author={Laiu, M. Paul and Tits, Andr{\'e} L.},
  journal={Computational Optimization and Applications},
  volume={72},
  pages={727--768},
  year={2019},
  publisher={Springer}
}

@article{schittkowski2009anactive,
title = {An active set strategy for solving optimization problems with up to 200,000,000 nonlinear constraints},
journal = {Applied Numerical Mathematics},
volume = {59},
pages = {2999--3007},
year = {2009},
author = {Klaus Schittkowski},
}

@article{facchinei1998onthe,
  title={On the accurate identification of active constraints},
  author={Facchinei, Francisco and Fischer, Andreas and Kanzow, Christian},
  journal={SIAM Journal on Optimization},
  volume={9},
  pages={14--32},
  year={1998},
  publisher={SIAM}
}

@article{duchi2021asymptotic,
 author = {John C. Duchi and Feng Ruan},
 journal = {The Annals of Statistics},
 pages = {21--48},
 publisher = {Institute of Mathematical Statistics},
 title = {Asymptotic optimality in stochastic optimization},
 volume = {49},
 year = {2021}
}

@incollection{polyak2001random,
title = {Random algorithms for solving convex inequalities},
publisher = {Elsevier},
volume = {8},
pages = {409--422},
year = {2001},
booktitle = {Inherently Parallel Algorithms in Feasibility and Optimization and their Applications},
author = {Boris T. Polyak},
}

@article{wei2020arandomized,
      title={A randomized nonlinear rescaling method in large-scale constrained convex optimization}, 
      author={Bo Wei and William B. Haskell and Sixiang Zhao},
      year={2020},
      journal={arXiv preprint arXiv:2003.10888},
}

@article{sonthalia2022project,
  author  = {Rishi Sonthalia and Anna C. Gilbert},
  title   = {Project and forget: {S}olving large-scale metric constrained problems},
  journal = {Journal of Machine Learning Research},
  year    = {2022},
  volume  = {23},
  pages   = {1--54},
}

@article{gould2005numerical,
title={Numerical methods for large-scale nonlinear optimization},
volume={14},
journal={Acta Numerica},
author={Gould, Nick and Orban, Dominique and Toint, Philippe},
year={2005},
pages={299--361}
}

@article{nedic2011random,
  title={Random algorithms for convex minimization problems},
  author={Nedi{\'c}, Angelia},
  journal={Mathematical Programming},
  volume={129},
  pages={225--253},
  year={2011},
  publisher={Springer}
}

@article{wardi1988stochastic,
  title={A stochastic algorithm for optimization problems with continua of inequalities},
  author={Wardi, Yorai},
  journal={Journal of Optimization Theory and Applications},
  volume={56},
  pages={285--311},
  year={1988},
  publisher={Springer}
}

@article{li2025new,
  title={New penalized stochastic gradient methods for linearly constrained strongly convex optimization},
  author={Li, Meng and Grigas, Paul and Atamt{\"u}rk, Alper},
  journal={Journal of Optimization Theory and Applications},
  volume={205},
  pages={29},
  year={2025},
  publisher={Springer}
}

@article{clarkson1995vegas,
  title={Las {V}egas algorithms for linear and integer programming when the dimension is small},
  author={Clarkson, Kenneth L.},
  journal={Journal of the ACM},
  volume={42},
  pages={488--499},
  year={1995},
  publisher={ACM}
}

@article{clarkson2012sublinear,
  title={Sublinear optimization for machine learning},
  author={Clarkson, Kenneth L. and Hazan, Elad and Woodruff, David P.},
  journal={Journal of the ACM},
  volume={59},
  pages={1--49},
  year={2012},
  publisher={ACM}
}

@article{burke1988on,
author = {Burke, James V. and Mor{\'e}, Jorge J.},
title = {On the identification of active constraints},
journal = {SIAM Journal on Numerical Analysis},
volume = {25},
pages = {1197--1211},
year = {1988},
}

@article{heunis1985use,
  title={Use of a {M}onte {C}arlo method in an algorithm which solves a set of functional inequalities},
  author={Andrew J. Heunis},
  journal={Journal of Optimization Theory and Applications},
  year={1985},
  volume={45},
  pages={89--99},
}

@inproceedings{clarkson1988las,
    author = {Clarkson, Kenneth L.},
    title = {A {L}as {V}egas algorithm for linear programming when the dimension is small},
    year = {1988},
    publisher = {IEEE},
    booktitle = {29{th} Annual Symposium on Foundations of Computer Science},
    pages = {452--456},
}

@article{tits2006constraint,
  title={Constraint reduction for linear programs with many inequality constraints},
  author={Tits, Andr{\'e} L. and Absil, Pierre-Antoine and Woessner, William P.},
  journal={SIAM Journal on Optimization},
  volume={17},
  pages={119--146},
  year={2006},
  publisher={SIAM}
}

@article{wang2015random,
  title={Random multi-constraint projection: {S}tochastic gradient methods for convex optimization with many constraints},
  author={Mengdi Wang and Yichen Chen and Jialin Liu and Yuantao Gu},
  journal = {arXiv preprint arXiv:1511.03760},
  year={2015},
}

@article{calafiore2010random,
  title={Random convex programs},
  author={Calafiore, Giuseppe Carlo},
  journal={SIAM Journal on Optimization},
  volume={20},
  pages={3427--3464},
  year={2010},
  publisher={SIAM}
}

@inproceedings{chakraborty2026randomized,
  title = 	 {Randomized Feasibility Methods for Constrained Optimization with Adaptive Step Sizes},
  author =       {Chakraborty, Abhishek and Nedi{\'c}, Angelia},
  booktitle = 	 {43{rd} International Conference on Machine Learning},
  pages = 	 {12475--12516},
  year = 	 {2026},
  volume = 	 {306},
  publisher =    {PMLR},
}

\appendix

\section{Hyperparameter Settings}

In this appendix, we provide the hyperparameters of the individual methods and Algorithm~\ref{fig:CST} used for the numerical experiments in Section~\ref{sec:empirical_results}, including the hyperparameters of the problem instances in Section~\ref{sec:dubins}.

For the experiments in Section~\ref{section:numerical_results}, for Algorithm~\ref{fig:CST}, we set~${\eta = 0.5}$, ${\tau_k = (k+1)^{0.6}}$, and $\varepsilon_k = 10^{-8}(k+1)^{-0.39}$ for all $k \in \N_0$, for every method and every problem instance. Moreover, for ABSGD \eqref{alg:absgd}, we set~${\delta_k^\mathrm{AB} = \hat{\delta}_k^\mathrm{AB} + \delta_\infty}$, where~$\{\hat{\delta}_k^\mathrm{AB}\}_{k\geq0}$ is a positive-valued, vanishing sequence as specified in Table~\ref{tab:numerical} below, and $\delta_\infty = 10^{-10}$ is a minimal positive floor added to avoid numerical artifacts. The hyperparameters for the methods~{\eqref{alg:projection}--\eqref{alg:apdsg}} are given in Table~\ref{tab:numerical}. Each method's hyperparameters are tuned to the best of our effort.

For the experiments in Section~\ref{sec:dubins}, we set $\eta =0.5$ and $\varepsilon_k = 10^{-8}(k+1)^{-0.1}$, and tune~${\omega >0}$ for~${\tau_k = (k+1)^{\omega}}$ for all $k \in \N_0$, for every method and every problem instance. Similar to the setting for Section~\ref{section:numerical_results}, for ABSGD \eqref{alg:absgd}, we set~${\delta_k^\mathrm{AB} = \hat{\delta}_k^\mathrm{AB} + \delta_\infty}$, where~$\{\hat{\delta}_k^\mathrm{AB}\}_{k\geq0}$ is a positive-valued, vanishing sequence as specified in Table~\ref{tab:dubins} below, and $\delta_\infty = 10^{-8}$ is a minimal positive floor added to avoid numerical artifacts. The hyperparameters given in Table~\ref{tab:dubins} are rounded to three significant figures. Each method's hyperparameters are tuned to the best of our effort. The hyperparameters of the problem instances for the experiments in Section~\ref{sec:dubins} are given in Table~\ref{tab:dubins_problems}.

\begin{table}[h]
\centering
\footnotesize
\caption{Algorithm hyperparameters for the experiments in Section~\ref{section:numerical_results}}
\label{tab:numerical}
\setlength{\tabcolsep}{3pt}
\begin{tabular}{@{}
    p{0.2\linewidth}
    p{0.15\linewidth}
    p{0.6\linewidth}
@{}}
\toprule
 Algorithm 
& Sampling 
& Hyperparameters \\
\midrule

ABSGD
& uniform
& \({\gamma_k^{\mathrm{AB}}}=\frac{1}{3}({k+1})^{-0.9}\), \(\hat{\delta}_k^\mathrm{AB}=\frac{1}{3}{({k+1})^{-0.9}}\) \\

& Algorithm~\ref{fig:CST}
& \({\gamma_k^{\mathrm{AB}}}={({k+1})^{-0.9}}\), \(\hat{\delta}_k^\mathrm{AB}={({k+1})^{-0.9}}\) \\

\midrule

Projection-SGD
& uniform
& \(\alpha_k=\frac{1}{3}{({k+1})^{-0.9}}\), ${\beta^{\Pi}_k} = 1$ \\

& Algorithm~\ref{fig:CST}
& \(\alpha_k={({k+1})^{-0.9}}\), ${\beta^{\Pi}_k} = 1$ \\

\midrule

Huber-penalty SGD
& uniform
& \(s_k=\frac{1}{3}{({k+1})^{-0.9}}\), \({\gamma_k^{\mathrm{H}}}=2\cdot 10^5 \log^2(k+1)\), \(\delta_k^\mathrm{H}={({k+1})^{-2}}\) \\

& Algorithm~\ref{fig:CST}
& \(s_k=({k+1})^{-0.9}\), \({\gamma_k^{\mathrm{H}}}=\frac{100}{3} \log^2(k+1)\), \(\delta_k^\mathrm{H}={({k+1})^{-2}}\) \\

\midrule

APDSG
& uniform
& \(\bar{\alpha}_k=3{({k+1})^{-0.9}}\), \(\varrho=10^4\), \(\bar{\rho}_k=\varrho\), \(\bar{\eta}=10^4\) \\

& Algorithm~\ref{fig:CST}
& \(\bar{\alpha}_k=3{({k+1})^{-0.9}}\), \(\varrho=10^3\), \(\bar{\rho}_k=\varrho\), \(\bar{\eta}=10^3\) \\

\bottomrule
\end{tabular}
\end{table}
\begin{table}[h]
\centering
\footnotesize
\caption{Algorithm hyperparameters for the experiments in Section~\ref{sec:dubins}}
\label{tab:dubins}
\setlength{\tabcolsep}{3pt}
\begin{tabular}{@{}
    p{0.17\linewidth}
    p{0.15\linewidth}
    p{0.11\linewidth}
    p{0.52\linewidth}
@{}}
\toprule
Problem & Algorithm
& Sampling 
& Hyperparameters \\
\midrule

Dubins Car 4.2~i) & ABSGD
& uniform
& ${\gamma_k^{\mathrm{AB}}}=0.0690{({k+1})^{-0.583}}$, $\hat{\delta}_k^\mathrm{AB}=0.0379{({k+1})^{-0.583}}$ \\

 & 
& Algorithm~\ref{fig:CST}
& ${\gamma_k^{\mathrm{AB}}}=0.110{({k+1})^{-0.560}}$, $\hat{\delta}_k^\mathrm{AB}=0.0657{({k+1})^{-0.560}}$,\\&&&$\tau_k = (k+1)^{0.178}$ \\

 & Projection-SGD
& uniform
& $\alpha_k=1.79{({k+1})^{-0.812}}$, ${\beta^{\Pi}_k} = 1$ \\

 & 
& Algorithm~\ref{fig:CST}
& $\alpha_k=0.793{({k+1})^{-0.676}}$, $ {\beta^{\Pi}_k} = 1$, $ \tau_k = (k+1)^{0.109}$ \\

\midrule

Dubins Car 4.2~ii) & ABSGD
& uniform
& ${\gamma_k^{\mathrm{AB}}}=0.0482{({k+1})^{-0.553}}$, $\hat{\delta}_k^\mathrm{AB}=0.0262{({k+1})^{-0.553}}$ \\

 & 
& Algorithm~\ref{fig:CST}
& ${\gamma_k^{\mathrm{AB}}}=0.0664{({k+1})^{-0.526}}$, $\hat{\delta}_k^\mathrm{AB}=0.0373{({k+1})^{-0.526}}$,\\
&&&$ \tau_k = (k+1)^{0.136}$ \\

 & Projection-SGD
& uniform
& $\alpha_k=1.77{({k+1})^{-0.821}}$, ${\beta^{\Pi}_k} = 1$ \\

 & 
& Algorithm~\ref{fig:CST}
& $\alpha_k=0.723{({k+1})^{-0.682}}$, ${\beta^{\Pi}_k} = 1$, $\tau_k = (k+1)^{0.0969}$ \\

\midrule

Dubins Car 4.2~iii) & ABSGD
& uniform
& ${\gamma_k^{\mathrm{AB}}}=0.000382{({k+1})^{-0.239}}$,\\&&&$\hat{\delta}_k^\mathrm{AB}=0.000201{({k+1})^{-0.239}}$ \\

 & 
& Algorithm~\ref{fig:CST}
& ${\gamma_k^{\mathrm{AB}}}=0.245{({k+1})^{-0.621}}$, $\hat{\delta}_k^\mathrm{AB}=0.149{({k+1})^{-0.621}}$,\\
&&&$\tau_k = (k+1)^{ 0.241}$ \\

 & Projection-SGD
& uniform
& $\alpha_k=1.78{({k+1})^{-0.815}}$, ${\beta^{\Pi}_k} = 1$ \\

 & 
& Algorithm~\ref{fig:CST}
& $\alpha_k=0.727{({k+1})^{-0.690}}$, ${\beta^{\Pi}_k} = 1$, $\tau_k = (k+1)^{ 0.109}$ \\

\midrule

Dubins Car 4.2~iv) & ABSGD
& uniform
& ${\gamma_k^{\mathrm{AB}}}=0.148{({k+1})^{-0.646}}$, $\hat{\delta}_k^\mathrm{AB}=0.0852{({k+1})^{-0.646}}$ \\

 & 
& Algorithm~\ref{fig:CST}
& ${\gamma_k^{\mathrm{AB}}}=0.0495{({k+1})^{-0.484}}$, $\hat{\delta}_k^\mathrm{AB}=0.0293{({k+1})^{-0.484}}$,\\
&&& $\tau_k = (k+1)^{0.144}$\\

 & Projection-SGD
& uniform
& $\alpha_k=1.78{({k+1})^{-0.806}}$, ${\beta^{\Pi}_k} = 1$ \\

 & 
& Algorithm~\ref{fig:CST}
& $\alpha_k=0.00326{({k+1})^{-0.309}}$, ${\beta^{\Pi}_k} = 1$, $\tau_k = (k+1)^{0.0982}$ \\

\bottomrule
\end{tabular}
\end{table}

\begin{table}[h]
\centering
\footnotesize
\caption{Problem hyperparameters for the experiments in Section~\ref{sec:dubins}}
\label{tab:dubins_problems}
\setlength{\tabcolsep}{3pt}
\begin{tabular}{@{}
    p{0.25\linewidth}
    p{0.7\linewidth}
@{}}
\toprule
Problem & Hyperparameters \\
\midrule

Dubins Car 4.2~i) & $N=3402$, $\epsilon_{\mathrm{safe}}=0.3$, $\epsilon_{\mathrm{unsafe}}=0.3$, $\epsilon_{\mathrm{dyn}}=0.1$, $b = 0.35$ \\

\midrule

Dubins Car 4.2~ii) & $N=3465$, $\epsilon_{\mathrm{safe}}=0.3$, $\epsilon_{\mathrm{unsafe}}=0.3$, $\epsilon_{\mathrm{dyn}}=0.1$, $b = 0.35$ \\

\midrule

Dubins Car 4.2~iii) & $N=3444$, $\epsilon_{\mathrm{safe}}=0.3$, $\epsilon_{\mathrm{unsafe}}=0.3$, $\epsilon_{\mathrm{dyn}}=0.1$, $b = 0.35$ \\

\midrule

Dubins Car 4.2~iv) & $N=3402$, $\epsilon_{\mathrm{safe}}=0.3$, $\epsilon_{\mathrm{unsafe}}=0.3$, $\epsilon_{\mathrm{dyn}}=0.1$, $b = 0.35$ \\

\bottomrule
\end{tabular}
\end{table}

\end{document}